\documentclass[12pt]{amsart}
\newtheorem{theorem}{Theorem}[section]
\newtheorem{proposition}[theorem]{Proposition}

\newtheorem{corollary}[theorem]{Corollary}
 
\theoremstyle{definition}

\newtheorem{remark}[theorem]{Remark}

\usepackage{amsfonts}
\usepackage{amssymb}
\usepackage{amsmath}
\usepackage{eucal}

\begin{document}
   
\begin{flushright}

[JAC2007(7RevisedLimit20)]\\
July 20, 2007

{\bf\it $''$ Never  shoot, Never  hit $''$}
  
\end{flushright}

\baselineskip=20pt

 \bigskip
\bigskip

\title{The Last Approach to the settlement of the Jacobian Conjecture}

\author{Susumu ODA}

\maketitle

\vspace{3mm}

\begin{abstract} The Jacobian Conjecture can be generalized and is established :  Let  $S$ be a polynomial ring over a field of characteristic zero in finitely may variables. Let  $T$  be an unramified, finitely generated  extension  of $S$ with  $T^\times = k^\times$.   Then  $T = S$.
\end{abstract}

\renewcommand{\normalbaselines}{\baselineskip20pt \lineskip3pt \lineskiplimit3pt}
\newcommand{\mapright}[1]{\smash{\mathop{\hbox to 1cm{\rightarrowfill}}\limits^{#1}}}

\newcommand{\lmapright}[1]{\smash{\mathop{\hbox to 7cm{\rightarrowfill}}\limits^{#1}}}

\newcommand{\mapleft}[1]{\smash{\mathop{\hbox to 1cm{\leftarrowfill}}\limits_{#1}}}

\newcommand{\lmapleft}[1]{\smash{\mathop{\hbox to 7cm{\leftarrowfill}}\limits_{#1}}}

\newcommand{\mapdown}[1]{\Big\downarrow\llap{$\vcenter{\hbox{$\scriptstyle#1\,$}}$ }}
\newcommand{\mapup}[1]{\Big\uparrow\rlap{$\vcenter{\hbox{$\scriptstyle#1$}}$ }}

\newcommand{\mapdownright}[1]{\searrow{\scriptstyle#1}}

\newcommand{\mapdownleft}[1]{{\scriptstyle#1}\swarrow}

\newcommand{\mapupright}[1]{\nearrow{\scriptstyle#1}}

\newcommand{\mapupleft}[1]{\nwarrow{\scriptstyle#1}}

\renewcommand{\thefootnote}{\fnsymbol{footnote}} %footnote


\vspace{4mm}  
              
 Let $ k $ be an algebraically closed field,
let $ \mathbb{A}^n_k $ be an affine space of dimension $ n $ over $ k $
and let  $ f : \mathbb{A}_k^n \longrightarrow \mathbb{A}_k^n $ be a morphism of affine spaces over $k$ of dimension $n$.  
Then $ f $  is given by coordinate functions $ f_1,  \ldots ,f_n,$
where $ f_i  \in k[X_1,  \ldots ,X_n] $ and
$ \mathbb{A}_k^n =  Max(k[X_1,\ldots,X_n]).$  If $ f $  has an inverse morphism,
then  the Jacobian $  \det( \partial f_i/ \partial X_j) $ is a nonzero
constant.  This follows from the easy chain rule.  The Jacobian Conjecture asserts the converse.  

\bigskip
       
If $ k $ is of characteristic \ $ p > 0 $  and $ f(X) = X  + X^p,$  then $ df/dX = f'(X) = 1 $ but
$ X $ can not be expressed as a polynomial in $ f$.
Thus we must assume the characteristic of $ k $ is zero.

\bigskip

\noindent
{\bf The Jacobian Conjecture of geometric form.}  \begin{sl}  Let  $f : \mathbb{A}_k^n \rightarrow \mathbb{A}_k^n$ 
 be a morphism of affine spaces of dimension $n\ (n \geq 1)$ over a field of characteristic zero.
 Then $f$ is expressed by coordinate functions  $ f_1,  \ldots ,f_n,$
where $ f_i  \in k[X_1,  \ldots ,X_n] $. 
 If  the Jacobian $  \det( \partial f_i/ \partial X_j) $ is a nonzero
constant,  then  $f$ is an isomorphism.
\end{sl}

\bigskip
The algebraic form of the Jacobian Conjecture is the following :
\bigskip

\noindent
{\bf The Jacobian Conjecture of algebraic form.}  \begin{sl}
 If  $f_1, \cdots, f_n$  be elements in a polynomial ring  $k[X_1, \cdots, X_n]$  over a field $k$ of characteristic zero such that  $  \det(\partial f_i/ \partial X_j) $ is a nonzero constant, then  $k[f_1, \cdots, f_n] = k[X_1, \cdots, X_n]$.
\end{sl}

\bigskip 

To prove the Jacobian Conjecture, we treat a more general case. More precisely, we show the following result:

\bigskip

\begin{sl}  
Let $k$ be a algebraically closed field of characteristic zero, let $S$ be a polynomial ring over $k$ of finite variables and let $T$ be an unramified, finitely generated extension domain of $S$ with $T^\times = k^\times$.    Then $T = S$.
\end{sl}

\bigskip
  
  Throughout this paper,  all fields, rings
and algebras are assumed to be  commutative with unity.
For a ring $ R,\ \   R^{\times} $ denotes the set of units of $ R $
and $ K(R) $ the total quotient ring.   ${\rm Spec}(R)$ denotes the affine scheme defined by $R$ or merely the set of all prime ideals of $R$ and ${\rm Ht}_1(R)$ denotes the set of all prime ideals of height one. 
Our general reference for unexplained technical terms is {\bf [9]}.

\vspace{30mm}
%%%
 
\section{Preliminaries} 

\bigskip
\bigskip

   \noindent
{\bf Definition.} Let $f : A \rightarrow B$ be a ring-homomorphism of finite type of locally Noetherian rings. The homomorphism $f$ is called {\it unramified} if $PB_P = (P\cap A)B_P$  and $k(P) = B_P/PB_P$ is a finite separable field extension of $k(P\cap A) = A_{P\cap A}/(P\cap A)A_{P \cap A}$  for all prime ideal $P$ of $B$.  The homomorphism $f$ is called {\it etale} if $f$ is unramified and flat.
 
 \bigskip

\begin{proposition}\label{2.1. Proposition}   Let  $k$  be an algebraically closed field of characteristic zero and let  $B$  be a polynomial ring  $k[Y_1, \ldots, Y_n]$.  Let  $L$  be a finite Galois extension of the quotient field of $B$  and  let  $D$  be an integral closure of  $B$  in $L$.  If  $D$  is etale over  $B$ then  $D = B$.
 \end{proposition}

 \begin{proof}  We may assume that  $k = {\bf C}$, the field of complex numbers  by "Lefschetz Principle"  (cf.{\bf [4, p.290]}).  The extension  $D/B$  is etale and finite, and so
 $$ {\rm Max}(D) \rightarrow {\rm Max}(B) \cong {\bf C}^n $$
 is a (connected) covering.  Since  ${\bf C}^n$  is simply connected, we have  $D = B$.  (An algebraic proof of the simple connectivity of $k^n$  is seen in {\bf [15]}.)  \end{proof}

\bigskip

Recall the following well-known results, which are required for proving Theorem \ref{2.2. Theorem} below.
\\
\\
{\bf Lemma A} ({\bf [9,(21.D)]}). \begin{it}    Let  $(A,m,k)$  and  $(B,n,k')$  be Noetherian local rings and $\phi : A \rightarrow B$  a local homomorphism $($i.e., $\phi(m) \subseteq n$ $)$.  If  $\dim B = \dim A + \dim B\otimes_Ak$  holds and if  $A$  and  $B\otimes_Ak = B/mB$  are regular, then  $B$  is flat over  $A$  and regular.    \end{it}

\bigskip

\noindent
{\it Proof.}  If  $\{\ x_1, \ldots, x_r\ \}$  is a regular system of parameters of  $A$ and if $y_1, \ldots, y
_s \in n$  are such that their images form a regular system of parameters of $B/mB$,  then $\{\ \varphi(x_1), 
\ldots, \varphi(x_r), y_1, \ldots, y_s\ \}$  generates $n$. and  $r+s = \dim B$. Hence $B$  is regular.
To show flatness, we have only to prove ${\rm Tor}^A_1(k,B) = 0$. The Koszul complex $K_*(x_1, \ldots, x_r;A)$
  is a free resolution of the $A$-module $k$. So we have ${\rm Tor}^A_1(k,B) = H_1(K_*(x_1, \ldots, x_r; A)\otimes_AB) = H_1(K_*(x_1, \ldots, x_r;B))$.  Since the sequence $\varphi(x_1), \ldots, \varphi(x_r)$ is a part of a regular system of parameters of $B$, it is a $B$-regular sequence.  Thus $H_i(K_*(x_1, \ldots, x_r;B)) = 0$
 for all $i > 0$. 
  \fbox

\bigskip

\noindent
{\bf Corollary A.1.} \begin{it} Let  $k$  be a field  and let $R = k[X_1, \ldots, X_n]$ be a polynomial ring.  Let  $S$  be a finitely generated ring-extension of  $R$.
If  $S$  is unramified over  $R$, then $S$  is etale over $R$.
  \end{it}

\bigskip

\noindent
{\it Proof.}  We have only to show that $S$  is flat over  $R$.  Take $P \in {\rm Spec}(S)$ and put $p=P \cap 
R$.  Then $R_p \hookrightarrow S_P$  is a local homomorphism. Since $S_P$  is unramified over  $R_p$,  we have
 $\dim S_P = \dim R_p$ and $S_P\otimes_{R_p}k(p) = S_P/PS_P = k(P)$ is a field. So by Lemma A, $S_P$  is flat over  $R_p$.  Therefore  $S$  is flat over $R$ by {\bf [5,p.91]}. 
 \fbox

\bigskip

\noindent
{\bf Example.}
  Let  $k$ be a field of characteristic $p > 0$ and let $S = k[X]$ be a polynomial ring.  Let $f = X + X^p \in S$.  Then the Jacobian matrix $\left(\dfrac{\partial f}{\partial X}\right)$ is invertible. So  $k[f] \hookrightarrow k[X]$ is finite and unramified.  Thus $k[f] \hookrightarrow k[X]$ is etale by Corollary A.1.  Indeed, it is easy to see that $k[X] = k[f] \oplus Xk[f] \oplus \cdots \oplus X^{p-1}k[f]$ as a $k[f]$-module, which implies that $k[X]$ is free over $k[f]$.

\bigskip

\noindent
{\bf Lemma B} ({\bf [2,Chap.V, Theorem 5.1]}). \begin{it} Let  $A$  be a Noetherian ring and  $B$  an $A$-algebra of finite type.  If  $B$  is flat over  $A$, then the canonical map  ${\rm Spec}(B) \rightarrow {\rm Spec}(A)$  is an open map.     \end{it}

\bigskip
 
 \noindent  
{\bf Lemma C} ({\bf [10, p.51,Theorem 3']}). \begin{it} Let  $k$  be a field and let $V$ be a $k$-affine variety defined by a $k$-affine ring $R$ (which means a finitely generated algebra over $k$) and let  $F$   be a closed subset of $V$ defined by an ideal $I$ of $R$.  If the variety $V \setminus F$  is $k$-affine, then $F$  is pure of codimension one.  \end{it}

\bigskip

%%%%%%%%%%%%%%%%
\noindent
{\bf Lemma D}({\bf [16,Theorem 9, $\S$ 4, Chap.V]}).  \begin{it}  Let $k$ be a field, let  $R$  be a $k$-affine domain and let  $L$ be a finite algebraic field extension of $K(R)$. Let $R_L$ denote the integral closure of $R$ in $L$. Then $R_L$ is a module finite type over  $R$. \end{it}

\bigskip
%%%%%%%%%%%%%

 %%%%%%%%%%%%%%%%%%

\noindent
{\bf Lemma E}({\bf [12, Ch.IV,Corollary 2]})(Zariski's Main Theorem).  \begin{it}  Let  $A$ be an integral domain and let $B$ be an $A$-algebra of finite type which is quasi-finite over  $A$.   Let  $\overline{A}$  be the integral closure of  $A$  in $B$.  Then  the canonical morphism  ${\rm Spec}(B) \rightarrow {\rm Spec}(\overline{A})$  is an open immersion.
  \end{it} 
 
 \bigskip

\noindent
{\bf Lemma F}({\bf [3, Corollary 7.10]}).  \begin{it} Let  $k$ be a field, $A$ a finitely generated $k$-algebra. Let  $M$  be a maximal ideal of  $A$. Then the field  $A/M$  is a finite algebraic extension of $k$.  In particular, if $k$ is algebraically closed then $A/M \cong k$.  \end{it}
 
\bigskip
%%%%%%%%%%%%%% 

\noindent
{\bf Lemma G }({\bf [2,VI(3.5)]}). \begin{it} Let $ f: A \rightarrow B$ and $g:B \rightarrow C$ be ring-homomorphisms of  finite type of locally Noetherian rings. 
 
{\rm (i)} Any immersion  ${}^af :{\rm Spec}(B) \rightarrow {\rm Spec}(A)$ is unramified.

{\rm (ii)} The composition $g\cdot f$ of unramified homomorphisms $f$ and $g$ is unramified.

{\rm (iii)} If $g\cdot f$ is an unramified homomorphism, then $g$ is an unramified homomorphism.
 \end{it}

\bigskip

\noindent
{\bf Lemma H }({\bf [2,VI(4.7)]}). \begin{it} Let $ f: A \rightarrow B$ and $g:B \rightarrow C$ be ring-homomorphisms of  finite type of locally Noetherian rings. $B$ (resp. $C$) is considered to be an $A$-algebra by $f$ (resp. $g\cdot f$).
 
{\rm (i)} The composition $g\cdot f$ of etale homomorphisms $f$ and $g$ is etale.
 
{\rm (ii)} Any base-extension $f\otimes_A1_C : C = A\otimes_AC \rightarrow B\otimes_AC$ of an etale homomorphism  $f$ is etale.

{\rm (iii)} If $g\cdot f : A \rightarrow B \rightarrow C$ is an etale homomorphism and if $f$ is an unramified homomorphism, then $g$ is etale.
 \end{it}

%%%%%%%%%%%%%%%%%%%%
\bigskip

\noindent
{\bf Corollary H.1.}  \begin{it} Let  $R$ be a ring and let $B \rightarrow C$ and $D \rightarrow E$ be etale $R$-algebra homomorphisms.  Then the homomorphism $B\otimes_RD \rightarrow C\otimes_RE$ is an etale homomorphism.
\end{it}

\bigskip

\noindent
{\it Proof.}  
The homomorphism
 $$ B\otimes_RD \rightarrow B\otimes_RE \rightarrow C\otimes_RE$$
 is given by composite of base-extensions.  So by Lemma H, this composite homomorphism is etale.
 \fbox

\bigskip
 
%%%%%%%%%%%%%  

\noindent
{\bf Lemma I }({\bf [11,(41.1)]})(Purity of branch loci).  \begin{it} Let  $R$  be a regular  ring and let $A$ be a normal ring which is a  finite extension of $R$. Assume that  $K(A)$ is finite separable extension of $K(R)$.  If $A_P$ is unramified over $R_{P\cap R}$  for all $P \in {\rm Ht}_1(A) (= \{ Q \in {\rm Spec}(A) | {\rm ht}(Q) = 1 \})$, then $A$ is unramified over $R$.  \end{it}

\bigskip

%%%%%%%%%%%%%%%%%%
\bigskip

\noindent
{\bf Lemma J} (cf. {\bf [17,(1.3.10)]}). \begin{it}  Let  $S$  be a scheme and let  $(X,f)$  and  $(Y,g)$  be $S$-schemes.  For a scheme  $Z$,  $|Z|$  denotes its underlying topological space.  Let  $p : X \times_S Y \rightarrow X$  and  $q : X \times_S Y \rightarrow Y$  be projections.  Then the map of topological spaces  $|p| \times_{|S|} |q| : |X \times_S Y| \rightarrow |X| \times_{|S|} |Y| $  is  a surjective map. 
   
\end{it}
 
\bigskip

\noindent
{\it Proof.}    Let $x \in X,\ y \in Y$ be points such that $f(x) = g(y) = s \in S$. Then the residue class fields $k(x)$ and $k(y)$ are the extension-fields of $k(s)$.  Let $K$  denote an extension-field of $k(s)$ containing two fields which are isomorphic to $k(x)$ and $k(y)$.  Such field $K$ is certainly exists.  For instance, we have only to consider the field $\mathcal{O}_{X,x}\otimes_{\mathcal{O}_{S,s)}}\mathcal{O}_{Y,y}/m$, where $m$  is a maximal ideal of  $\mathcal{O}_{X,x}\otimes_{\mathcal{O}_{S,s)}}\mathcal{O}_{Y,y}$.  Let $x_K : {\rm Spec}(K) \rightarrow  {\rm Spec}(\mathcal{O}_{X,x}) \mapright{i_x} X$, where  $i_x $ is the canonical immersion as  topological spaces and the identity  $i_x^*(\mathcal{O}_X) = \mathcal{O}_{X,x}$ as structure sheaves. Let $y_K$  be the one similarly defined as $x_K$.  By the construction of $x_K,\ y_K$,  we have $f\cdot x_K = g\cdot y_K$. Thus there exists a $S$-morphism $z_K: {\rm Spec}(K) \rightarrow X\times_SY$  such that $p\cdot z_K = x_K,\ q\cdot z_K = y_K$.  Since ${\rm Spec}(K)$ consists of a single point, putting its image $= z$, we have $p(z) = x,\ q(z) = y$.  Therefore the map of topological spaces  $|p| \times_{|S|} |q| : |X \times_S Y| \rightarrow |X| \times_{|S|} |Y| $  is surjective.

    \fbox

%%%%%%%%%%%%%%%%%%
\bigskip

\noindent
{\bf Remark 1.1.}  Let $A \rightarrow B$ be a ring-homomorphism of rings.  Let $pr_i : {\rm Spec}(B) \times_{{\rm Spec}(A)}{\rm Spec}(B) \rightarrow {\rm Spec}(B)\ (i=1,2)$ be the projection.
Recall that an affine scheme ${\rm Spec}(B)$ is separated over ${\rm Spec}(A)$, that is, the diagonal morphism $\Delta : {\rm Spec}(B) \rightarrow {\rm Spec}(B) \times_{{\rm Spec}(A)}{\rm Spec}(B)$ (defined by $B\otimes_AB \ni x\otimes y \mapsto xy \in B$) is a closed immersion and  $pr_i\cdot \Delta = id_{{\rm Spec}(B)}\ (i=1.2)$ \ (cf. {\bf [17]}).   
  It is easy to see that the diagonal morphism $\Delta' :  {\rm Spec}(B) \rightarrow {\rm Spec}(B) \times_{{\rm Spec}(A)} \cdots \times_{{\rm Spec}(A)}{\rm Spec}(B)\ (n\mbox{-times})$ similarly defined is also a closed immersion with $p_i\cdot \Delta' = id_{{\rm Spec}(B)}$, where $pr_i$ is the projection\ $(1 \leq i \leq n)$.
  Let $B_2, \ldots, B_n$ be $A$-algebras such that  $B \cong_A B_2 \cong_A \cdots  \cong_A B_n$.
 Then there exists a ${\rm Spec}(A)$-morphism $\Delta^* : {\rm Spec}(B) \rightarrow {\rm Spec}(B) \times_{{\rm Spec}(A)} \cdots \times_{{\rm Spec}(A)}{\rm Spec}(B) \cong_{{\rm Spec}(A)} {\rm Spec}(B) \times_{{\rm Spec}(A)} {\rm Spec}(B_2) \times_{{\rm Spec}(A)} \cdots  \times_{{\rm Spec}(A)} {\rm Spec}(B_n)$, which is a closed immersion and $pr_1\cdot \Delta^* = id_{{\rm Spec}(B)}$. Hence  $pr_1$ is surjective.

\bigskip

\noindent
{\bf Remark 1.2.} Let $k$ be a field, let  $ S = k[Y_1,\ldots,Y_n]$ be a polynomial ring over $k$ and let $L$ be a finite Galois extension field of $K(S)$ with Galois group $G = \{\ \sigma_1 = 1, \sigma_2, \ldots, \sigma_\ell\}$.
  Let  $T$ be a finitely generated, flat extension of $S$ contained in $L$ with $T^\times = k^\times$. Put $T^{\sigma_i} = \sigma_i(T) \subseteq L$.
  Let 
 $$T^\# := T^{\sigma_1} \otimes_S \cdots \otimes_S T^{\sigma_\ell},$$
 which has the natural $T$-algebra structure  by $T\otimes_SS\otimes_S \cdots \otimes_SS \hookrightarrow T^{\sigma_1} \otimes_S \cdots \otimes_S T^{\sigma_\ell} = T^\#$.  

(i) Let $P$ be a prime ideal of $T$.
  Then the element  $(P^{\sigma_1}, \ldots, P^{\sigma_\ell}) \in |{\rm Spec}(T^{\sigma_1})|$  $ \times_{|{\rm Spec}(S)|} \cdots \times_{|{\rm Spec}(S)|}|{\rm Spec}(T^{\sigma_\ell})|$  is an image of some element $Q$ in $|{\rm Spec}(T^\#)|$  because   the canonical map $ |{\rm Spec}(T^\#)| = |{\rm Spec}(T^{\sigma_1} \otimes_S \cdots \otimes_S T^{\sigma_\ell})|  \rightarrow  |{\rm Spec}(T^{\sigma_1})| \times_{|
{\rm Spec}(S)|} \cdots \times_{|{\rm Spec}(S)|}|{\rm Spec}(T^{\sigma_\ell})|$ is surjective by Lemma J. The map $|{\rm Spec}(T^\#)| \rightarrow |{\rm Spec}(T)|$  yields that $Q \cap T = P$  
 Hence $|{\rm Spec}(T^\#)| \rightarrow |{\rm Spec}(T)|$ is surjective. (This result has been obtained in Remark 1.1.) 
  So  $T^\#$ is faithfully flat over $T$.

 (ii) Take  $p \in {\rm Ht}_1(S)$.  Then $p$ is a principal ideal of $S$ and so $pT^{\sigma_i} \not= T^{\sigma_i}\ (\forall \sigma_i \in G)$  because $T^\times = k^\times$.
   Let $P$ be a minimal prime divisor of $pT$.   Then $P^{\sigma_i} \in {\rm Spec}(T^{\sigma_i})$ and  $P^{\sigma_i} \cap S = p$ because $S \hookrightarrow T$ is flat.
  There exists a prime ideal $Q$ in ${\rm Spec}(T^\#)$ with $Q \cap T = P$ by (i)   and hence  $P \cap S = p$.  Thus $Q \cap S = p$. 
   Therefore $pT^\# \not= T^\#$ for all $p \in {\rm Ht}_1(S)$.  

\bigskip

%%%%%%%%%%%  
%%%%%%%%%%   

\vspace{20mm}
%%%
 
\section{Main Result} 
 
 \bigskip 

The following is our main theorem.

\begin{theorem}\label{2.2. Theorem}
 Let $k$ be a algebraically closed field of characteristic zero, let $S$ be a polynomial ring over $k$ of finitely many variables and let $T$ be an unramified, finitely generated extension domain of $S$ with $T^\times = k^\times$.   Then $T = S$.
%%%%%%%%%%%
\end{theorem}

\noindent
{\it Proof.}   

{\bf (1)}  Let  $K( \ \ )$  denote the quotient field of  $(\ \ )$.  There exists a minimal finite Galois extension  $L$ of  $K(S)$  containing  $T$  because  $K(T)/K(S)$  is a finite algebraic extension.
  
  Let  $G$ be the Galois group  $G(L/K(S))$.  
Put  $G = \{\ \sigma_1 = 1, \sigma_2, \ldots, \sigma_\ell\}$,  where  $\sigma_i \not=
 \sigma_j$  if  $i \not= j$.    Put $T^{\sigma} := \sigma(T) \ (\forall \ \sigma \in 
G)$  and put $D := S[\bigcup_{\sigma \in G}T^{\sigma}] = S[\bigcup_{i=1}^\ell T^{\sigma_i}] \subseteq L$.  Then  $K(D) = L$ since $L$ is  a minimal Galois extension of $K(S)$ containing $K(T)$.   
 Since  ${\rm Spec}(T) \rightarrow {\rm Spec}(S)$  is etale (Corollary A.1 or {\bf [4, p.296]}),  so is  ${\rm Spec}(T^{\sigma}) \rightarrow {\rm Spec}(S)$  for each  $\sigma \in G$.
     
 Put 
 $$T^\# := T^{\sigma_1} \otimes_S \cdots \otimes_S T^{\sigma_\ell},$$
 which has the natural $T$-algebra structure  by $T = T\otimes_SS\otimes_S \cdots \otimes_SS \hookrightarrow T^{\sigma_1} \otimes_S \cdots \otimes_S T^{\sigma_\ell} = T^\#$.
  This homomorphism is etale by Corollary H.1 because $S \rightarrow T$ is etale.
  Let  $\psi' :  T^\# = T^{\sigma_1} \otimes_S \cdots \otimes_S T^{\sigma_\ell} \rightarrow L$  be an $S$-algebra homomorphism sending  $a_1^{\sigma_1}\otimes \cdots \otimes a_\ell^{\sigma_\ell}$  to  $a_1^{\sigma_1}\cdots a_\ell^{\sigma_\ell}\ (a_i \in T)$.  Then  $D = {\rm Im}(\psi') = S[\bigcup_{\sigma \in G}T^{\sigma}] \subseteq L$.  
%%% 
  %%%
  Since  ${\rm Spec}(T) \rightarrow {\rm Spec}(S)$  is  etale,  the canonical morphism   $ {\rm Spec}(T^\#) = {\rm Spec}(T^{\sigma_1} \otimes_S \cdots \otimes_S T^{\sigma_\ell})  \rightarrow  {\rm Spec}(T^{\sigma_1}\otimes_SS\otimes_S \cdots \otimes_SS) = {\rm Spec}(T)$ is etale, and the natural surjection  $\psi : T^\# = T^{\sigma_1} \otimes_S \cdots \otimes_S T^{\sigma_\ell} \rightarrow D$  is unramified by Lemma G(i)(or {\bf [2,VI(3.5)]}).
 So $[T \hookrightarrow D] = [T \hookrightarrow T^\# \rightarrow D]$  is unramified by Lemma G(ii) because etale is flat and unramified. Moreover $S \hookrightarrow T \hookrightarrow D$ is also unramified.
    Since  $T$ and $D$ are unramified over  $S$, both $T$ and $D$  are etale over  $S$  and  both $T$ and $D$ are regular by Corollary A.1.

 Let  $I := \text{Ker} \psi$.  So $ {}^a\psi : {\rm Spec}(D) \cong V(I) \subseteq {\rm Spec}(T^\#)$  is a closed immersion.  Since  $[T \hookrightarrow T^\# \rightarrow D] = [T \hookrightarrow D]$  is etale,  so is  $ \psi : T^\# \rightarrow D$  by Lemma H(iii) (or {\bf [2,VI(4.7)]}).
     It follows that  ${\rm Spec}(D) \rightarrow {\rm Spec}(T^\#)$  is a closed immersion and an open map because  a flat morphism is an open map by Lemma B.   Thus  ${\rm Spec}(D) = V(I) \subseteq {\rm Spec}(T^\#)$  is a connected component of ${\rm Spec}(T^\#)$.
 So we have seen  that the natural $S$-homomorphism  $T \hookrightarrow T^{\#} \rightarrow D$  is etale and that  ${\rm Spec}(D)$  is a connected component of ${\rm Spec}(T^{\#})$. 
   Note that  $T^\#$  is reduced  because  $T^\#$ is  unramified over  $S$, and that   $\dim S = \dim T = \dim D$  because $S, T$ and $D$ are all $k$-affine domains with the same transcendence degree over  $k$.
 
   Let  $(0) = \bigcap_{i=1}^s P_i $  be an irredundant primary decomposition in $T^\#$.
   %%%%%%% 
    Since $T \rightarrow T^\#$  is  flat,  the GD-theorem {\bf [9,(5.D)]}(or Lemma B) holds for this homomorphism  $T \rightarrow T^\#$. In the decomposition $(0) = \bigcap_{i=1}^sP_i$, each $P_i$  is a minimal prime divisor of $(0)$, so we have  $T \cap P_i = (0)$ for all $i=1, \ldots,s$.   Note that  $S \hookrightarrow T^{\sigma_i}$ is unramified and hence that  $T^\#$  is reduced. The $P_i$'s  are  prime ideals of  $T^\#$. 
      Note that $I$  is a prime ideal of  $T^\#$  and that  $\dim S = \dim T = \dim T^{\sigma} = \dim D$   for each $\sigma \in G$.  
 Thus  there exists  $j$, say $j = 1$,  such that  $I = P_1$.
 In this case, $P_1 + \bigcap_{i=2}^s P_i = T^\#$  and  $T^\#/P_1 \cong D \subseteq L$ as $T$-algebra. 
  Note that  $T$  is considered to be a subring of $T^\#$  by the canonical injective  homomorphisms $ T = T \otimes_SS\otimes_S\cdots \otimes_SS \hookrightarrow T^\#$ and that $[T \hookrightarrow T^\# \rightarrow T^\#/P_1 \cong D] = [T \hookrightarrow D]$.
 Putting $C = T^\#/\bigcap_{i=2}^sP_i$,  we have 
  $  T^\# \stackrel{\Phi}{\widetilde{\rightarrow}} T^\#/P_1 \times T^\#/\bigcap_{i=2}
^sP_i \cong D \times C.$  
 The ring  $D$ is considered  a $T$-algebra naturally and $D \cong_T T^\#/P_1$. 
  Similarly we can see that $P_i + P_j = T^\#$ for any $i \not= j$. So  consider $T^\#/P_j$ instead of $D$, we have a direct product decomposition:
  $$ \Phi : T^\# \cong T^\#/P_1 \times \cdots \times T^\#/P_s.$$  
    Considering $ T = T\otimes_SS\otimes_S \cdots \otimes_SS \hookrightarrow T^{\sigma_1} \otimes_S \cdots \otimes_S T^{\sigma_\ell} = T^\# \rightarrow T^\#/P_i\ (1 \leq i \leq s)$,  $T^\#/P_1$ is a $T$-algebra $(1 \leq i \leq s)$ and  $\Phi$ is a $T$-algebra isomorphism.  Moreover  each $T^\#/P_i$ is regular (and hence normal) and no non-zero element of $T$ is  a zero-divisor on $T^\#/P_i\ (1\leq i \leq s)$.

\bigskip
    
%%%%%%%%%%%%%%%%%%%%%%%%%%%%%%%%  
{\bf (2)}  Now we claim that 
$$ aD \not= D\  \ (\forall a \in S \setminus S^\times) \ \ \ \ \ \ {\bf (\#)}.$$
 
 Note first that for all $p \in {\rm Ht}_1(S)$,  $pT \not= T$  because $p$ is principal and $T^\times = k^\times$,   and hence that $pT^\sigma \not= T^\sigma $  for all $\sigma \in G$.  Thus  $pT^\# \not= T^\#$  for all $p \in {\rm Ht}_1(S)$ by Remark 1.2.  Since $S$ is a polynomial ring, any $p \in {\rm Ht}_1(S)$ is principal.

   Let  $a \in S\ (\subseteq T^\#)$ be any non-zero prime element in $S$.  Then by    the above argument, $aT^\# \not= T^\#$.  When $s = 1$, then the assertion $(\#)$  holds  obviously.  So we may assume that $s \geq 2$.

 %%%%%%%%%
 
 {\bf  Suppose that $a \in S$ is a prime element and that  $aD = D$.}
   
 %%%%%%%%%%  
Then $aT^\# + P_1 = T^\#$ and $P_2\cdots P_s = T^\#(P_2\cdots P_s) = (aT^\# + P_1)(P_2 \cdots P_s) = aP_2\cdots P_s$  because  $P_1\cdots P_s = (0)$.  That is,

$$  aP_2\cdots P_s = P_2\cdots P_s\ \ \ \ \  \ \ \ \ \ \ (*).$$

  Throughout  this proof, for a subset $V$ of $T^\boxtimes$,  $V^\times$ denotes $T^\boxtimes \cap V$,

%%%%%%%%%%%%%%%%%%%%%%%%%%%%%%%%%%%%%%%%%%%%%%%%%%%%%%%%%%%%%%%%%%%%%
 %%%%%%%%%%%%%%%%%%%%%%%%%%%%%%%%%%%%%%%%%%%%%%%%%%%%%%%%%%%%%%%%%%%%

Put $ p = aS \in {\rm Ht}_1(S)$.  
  Let $T^\#_p := T^\# \otimes_SS_p = T_p^{\sigma_1} \otimes_{S_p} \cdots \otimes_{S_p}T_p^{\sigma_\ell}$, (which is a semi-local ring  because $S \rightarrow T^\#$ is etale).
 Note that the Going Up Theorem holds for $S_p \subseteq T_p$ because both $S$ and $T$ are integral domain and ${\rm ht}(p) = 1$. 
  Since $pT^\# \not= T^\#$, we have $pT^\#_p \not= T^\#_p$.

  Any prime ideal $P$ of $T^\#_p = (S\setminus p)^{-1}T^\#$ is $(P \cap T^\#)(S\setminus p)^{-1}T^\#$, that is, there exists the canonical bijection ${\rm Spec}((S\setminus p)^{-1}T^\#) \cong \{ Q \in {\rm Spec}(T^\#) | (S\setminus p) \cap Q = \emptyset \}$ corresponding $P \mapsto P \cap T^\#$.
  
  Let $M$ be a maximal ideal of $T^\#_p$.  Then $M' = M \cap T^\#$ is a prime ideal satisfying $M' \cap (S\setminus p) = \emptyset$.  So $M \cap S$ is either $(0)$ or $p$.  

Suppose that  $M \cap S = (0)$, that is, $M' \cap S = (0)$.  Then $M' \cap T = (0)$ and  ${\rm ht}(M') = 0$ because $T$ is algebraic over $S$ and  $S \rightarrow T^\#$ is etale.  
   Let $T^\boxtimes = T \otimes_S \cdots \otimes_ST\ (\ell\mbox{-times})$ and  $\lambda : T^\boxtimes \rightarrow T$ be an $S$-algebra homomorphism sending $c_1\otimes \cdots \otimes c_\ell$ to $c_1 \cdots c_\ell$  with $c_i \in T$.  The $S$-algebra  $T^\boxtimes$ can be $T$-algebra by the canonical homomorphism $T = T \otimes_S S \otimes_S \cdots \otimes_SS \rightarrow T^\boxtimes$.   
  Let $\Psi : T^\# \rightarrow T^\boxtimes$ be an $S$-isomorphism sending $c_1^{\sigma_1}\otimes \cdots \otimes c_\ell^{\sigma_\ell}$ to $c_1\otimes \cdots \otimes c_\ell$  with $c_i \in T$  and let $\Psi_p : T^\#_p \cong T^\boxtimes_p$.    Then $M''_p = \Psi(M'_p) = \Psi(M')_p$.   Note that $\lambda$ is an etale surjection.  
  Put $M'' = \Psi(M')$.
 Then $M'' \cap T = (0)$ in $T^\boxtimes$.  
  It is easy to see that the $S$-algebra homomorphisms $\Psi$ and $\lambda$ can be  $T$-algebra  homomorphisms in the natural way. 
 
(i) If $\lambda(M'') = T$, then the restriction $\lambda| : M'' \rightarrow T$ is a split surjection as $T$-modules.
  Then $T^\boxtimes/{\rm Ker}(\lambda) \cong_T \lambda(M'') = T$. 
 So $T^\boxtimes \cong_T{\rm Ker}(\lambda) + T$. 
 Thus $T^\boxtimes = {\rm Ker}(\lambda) + mT$ for some $m \in M''$ with $\lambda(m) =1$.  Since  for any  $t \in T$, $\lambda(mt - t) = \lambda(m)\lambda(t) - t = t - t = 0$,  we have $mT + {\rm Ker}(\lambda) = T + {\rm Ker}(\lambda)$. 
  Note that  both   ${\rm Ker}(\lambda)$ and $M''$ are contained in $\{ \Psi(P_1) ,\ldots, \Psi(P_s) \}$ since ${\rm ht}({\rm Ker}(\lambda)) = 0 = {\rm ht}(M'')$.
 So $M'' = \Psi(P_i)$ and ${\rm Ker}(\lambda) = \Psi(P_j)$.
   Since  $\lambda(M'') = T$, we may assume that $\Psi(P_2) = {\rm Ker}(\lambda)$, otherwise $D \cong_T  T^\boxtimes/{\rm Ker}(\lambda) \cong_T T$ and $D^\times = T^\times = k^\times$, a contradiction.
  Let  $\Psi' : T^\boxtimes \cong_T T^\boxtimes/M'' \times  T^\boxtimes/{\rm Ker}(\lambda) \times T^\boxtimes/\Psi(P_3) \times \cdots \times T^\boxtimes/\Psi(P_s) $ be the isomorphism induced from  $\Psi$. 
 %%%%%%%%%%%%%%%%%%%%%%%%%%
 Then $ \lambda(M'') = T = \lambda {\Psi'}^{-1}(T^\boxtimes/{\rm Ker}(\lambda)) = \lambda(T^\boxtimes) = \lambda ({\Psi'}^{-1}((T^\boxtimes/M'') \times  (T^\boxtimes/{\rm Ker}(\lambda)) \times (T^\boxtimes/\Psi(P_3)) \times \cdots \times (T^\boxtimes/\Psi(P_s)))$.  Hence we have  ${\Psi'}^{-1}((T^\boxtimes/\Psi(M'') \times 0 \times (T^\boxtimes/\Psi(P_3) \times \cdots \times T^\boxtimes/\Psi(P_s)) \subseteq {\rm Ker}(\lambda)$, but since $\lambda(M'')  = T$   and ${\rm Ker}(\lambda)$ is a prime ideal of  $T^\boxtimes$,  $\Psi^{-1}(T^\boxtimes/M'')$ must be ${\rm Ker}(\lambda)$.  Thus $T^\boxtimes/\Psi(P_3) \times \cdots \times T^\boxtimes/\Psi(P_s) = 0$, which means that $s = 2$ in this case.
%%%%%%%%%%%%%%%%%%%%%%
 So $M'' \cap \Psi(P_2) = (0)$.   Moreover  $T^\boxtimes = {\rm Ker}(\lambda) + mT = {\rm Ker}(\lambda) + M''$, it follows that $T^\boxtimes \cong_T T^\boxtimes/M'' \times T^\boxtimes/{\rm Ker}(\lambda) =  T^\boxtimes/M'' \times T^\boxtimes/\Phi(P_2)$.
  Thus $D \cong_TT^\boxtimes/M''$. 
 Thus 
     $$\Psi' : T^\boxtimes \cong_T T^\boxtimes/M'' \times  T^\boxtimes/{\rm Ker}(\lambda) \ \ \ \ \ \ \ \ \ \  \ \ \ \ \ \ (**)$$
 whence  $M'' = \Psi(P_1)$ and ${\rm Ker}(\lambda) = \Psi(P_2)$.  
 Here in this case,  $s = 2$  in (1).
    We have 
 $(T^\boxtimes)^\times \subseteq mT^\times + {\rm Ker}(\lambda) = mk^\times + {\rm Ker}(\lambda)= k^\times + {\rm Ker}(\lambda)$.
 Note that  $ (k^\times + {\rm Ker}(\lambda))^\times$ is a group by the multiplication in $T^\boxtimes$.
Thus 
   $$ (T^\boxtimes)^\times \subseteq (k^\times + {\rm Ker}(\lambda))^\times\ \ \ \ \ \ \ \ \ \ \  \ \ \ \ \ \ \ \ \ \ \ \ \  (***)$$
  From  $(**)$, we have  
$$k^\times \times k^\times \subseteq  (T^\boxtimes)^\times  \cong (T^\boxtimes/M'')^\times \times (T^\boxtimes/{\rm Ker}(\lambda))^\times\ \ \ \ \ \ \ \ \ \ \ (****)$$
%%%%%%%%%%%%%
 It is easy to see that  $(T^\boxtimes/{\rm Ker}(\lambda))^\times \cong T^\times  = k^\times$,
  and so that  by $(****)$  and $(***)$  we have  $(T^\boxtimes/M'')^\times \subseteq  (k^\times + {\rm Ker}(\lambda))^\times/k^\times  \subseteq (1 + {\rm Ker}(\lambda))^\times$, that is $k^\times \subseteq (1 + {\rm Ker}(\lambda))^\times$,  which is impossible  because if we take any $c \in k^\times$ then $1-c \in {\rm Ker}(\lambda) \cap k$ implies $c = 1$. 
%%%%%%%%%%%%%%%%%%%%% 
           So this case does not occur.

(ii) If $\lambda(M'') \cap S  = p$, then it is easy to see that $M' \cap S = p$, a contradiction.

(iii) Let $\lambda(M'') \cap S = (0)$. 
  In this case, $\lambda(M'')  = (0)$ because $S \hookrightarrow T$ is algebraic,  and hence $M'' \subseteq {\rm Ker}(\lambda)$.
   We have an $S_p$-isomorphism $T^\boxtimes_p/{\rm Ker}(\lambda)_p \cong T_p$.  Since  $pT_p \not= T_p$, there exists a prime ideal $N''$ of $T^\boxtimes$  such that $N'' \supset M''$, \ ${\rm ht}(N'') = 1$ and $N'' \cap S = p$  because $\lambda$ is etale. 
 So $N' := \Psi^{-1}(\lambda^{-1}(N''))$ satisfies $N'_p \supsetneq M'_p = M$ and $N'_p \cap S = p$ because $\lambda$ is etale, which  contradicts the maximality of $M$.
 
  Therefore $M' \cap S = M \cap S = p$.

%%%%%%%%%%%%%%%%%%%%%%%%%%%%%%%%%%%%%%%%%%%%%%%%%%%%%%%%%%%%%%%%%%%
%%%%%%%%%%%%%%%%%%%%%%%%%%%%%%%%%%%%%%%%%%%%%%%%%%%%%%%%%%%%%%%%%%%
  So we conclude that the Jacobson radical $J(T^\#_p)$ of $T^\#_p$ is $\sqrt{pT^\#_p}$ and  contains the prime element  $a$.

%%%%%%

From $(*)$, we have $a{P_2}_p \cdots {P_s}_p = {P_2}_p \cdots {P_s}_p$, which is a finitely generated $T^\#_p$-module.  
Thus there exists $\beta \in T^\#_p$  such that $(1 - a\beta){P_2}_p\cdots {P_s}_p= 0$.  Since $a$ is contained in the Jacobson radical $J(T^\#_p)$ of the semi-local ring $T^\#_p$ as mentioned above,  we have ${P_2}_p \cdots {P_s}_p = 0$.
 Since any element of $S \setminus p$ is not a zero-divisor on $T^\#$, we have $P_2 \cdots P_s \subseteq {P_2}_p \cdots {P_s}_p = (0)$.  So $ P_2 \cap \cdots \cap P_s = P_2 \cdots P_s = (0)$.
But $(0) = P_1\cap \cdots \cap P_s$ is an irredundant primary decomposition  as mentioned above,  which is a {\bf contradiction}.
      Hence $(\#)$ has been proved.

 \bigskip 
  %%%%%%% 
     %%%%%%%%%%%%%%%%%%%%%%%%%%%%%%%%% 
     
  {\bf (3)}     Let $C$  be the integral closure of $S$ in $L$.  Then $C \subseteq D$  because  $D$ is regular (hence normal) and $C$ is an $k$-affine domain (Lemma D). For any $\sigma \in G = G(L/K(S))$, $C^\sigma \subseteq D$  because  $C^\sigma$  is integral over $S$ and $D$ is normal with $K(C) = L$.  Hence   $C^\sigma = C$  for any  $\sigma \in G$. Note that both $D$  and $C$ have the quotient field $L$.  
   Zarisiki's  Main Theorem(Lemma C) yields the decomposition: 
$${\rm Spec}(D) \mapright{i} {\rm Spec}(C) \mapright{\pi} {\rm Spec}(S),$$
 where $i$ is an open immersion and $\pi$ is integral(finite).  
  We identify  ${\rm Spec}(D) \hookrightarrow {\rm Spec}(C) $ as open subset and $D_P = C_{P \cap C}\ ( P \in {\rm Spec}(D))$.  
 Let $Q \in {\rm Ht}_1(C)$ with $Q \cap S = p = aS$. Then $Q$ is a prime divisor of $aC$. Since  $aD \not= D$ by $(\#)$ in (2), there exists $P \in {\rm Ht}_1(D)$  such that $P \cap S = p$.  Hence   there exists $\sigma \in G$  such that $Q = (P \cap C)^\sigma$ because  any minimal divisor of $aC$ is $(P \cap C)^{\sigma'}$ for some $\sigma' \in G$ ({\bf [9,(5.E)]}), noting that $C$ is a Galois extension of $S$.  Since  $D_P = C_{P \cap C}$  is unramified over $S_{P \cap S} = S_p$, $C_Q = C_{(P\cap C)^\sigma} \cong C_{(P\cap C)}$  is unramified over $S_p$.  
 Hence $C$ is unramified over $S$ by Lemma I.  By Corollary A.1,  $C$ is  finite etale over $S$.  So Proposition \ref{2.1. Proposition}  implies that $C = S$.  In particular, $L = K(D) = K(C) = K(S)$ and hence $K(T) = K(S)$.  Since $S \hookrightarrow T$ is birational etale, ${\rm Spec}(T) \hookrightarrow {\rm Spec}(S)$ is an open immersion by Lemma C.  Let $J$  be an ideal of $S$ such that $V(J) = {\rm Spec}(S) \setminus {\rm Spec}(T)$.
  Suppose that $J \not= S$. 
 Then $V(J)$ is  pure of  codimension one by Lemma C. Hence $J $ is a principal ideal $aS$  because $S$ is a UFD.  Since $JT = aT = T$,  $a$ is a unit in $T$. 
 But $T^\times = k^\times $ implies that $a \in k^\times$ and hence that $J = S$, a contradiction.   Hence  $V(J) = \emptyset$, that is,  $T = S$.
   {\bf Q.E.D.}

 %%%%%%%%%%%%%%
 %%%%%%%%%%%%%%
 %%%%%%%%%%%%%%

 \vspace{15mm}

\section{The Jacobian Conjecture}

 The Jacobian conjecture has been settled affirmatively
in several  cases.  For example,\\
Case(1) $ k(X_1,\ldots ,X_n) $ is a Galois
extension of $ k(f_1,\ldots,f_n) $  (cf.{ \bf [4],[6]} and
{\bf [15]});\\
Case(2)  $ \deg f_i \leq 2 $\  for all $ i $ (cf.{ \bf [13]} and { \bf [14]});\\
Case(3) $ k[X_1,\ldots ,X_n] $ is integral over $ k[f_1,\ldots,f_n].$  (cf.{ \bf [4]}).

  A general reference for the Jacobian Conjecture is { \bf [4]}.

 \begin{remark}\label{2.4. Remark} 
 
 (1)  In order to prove Theorem \ref{2.3. Theorem}, we have only to show that the  inclusion $ k[f_1,\ldots,f_n]$ $\longrightarrow $ $k[X_1,\ldots,X_n]$  is surjective.  For this it suffices
 that $ k'[f_1,\ldots,f_n] \longrightarrow$ $k'[X_1,\ldots,X_n]$  is surjective, where  $k'$  denotes an algebraic closure of $ k$.  Indeed,  once we proved  $k'[f_1,\ldots,f_n]$ $ = k'[X_1,\ldots,X_n]$,  we can  write for each  $i = 1,\ldots,n$: 
  $$         X_i = F_i(f_1,\ldots,f_n),$$ 
 where $F_i(Y_1,\ldots,Y_n) \in k'[Y_1,\ldots,Y_n]$, a polynomial ring in $Y_i$. Let $ L$  be an intermediate field between $ k$  and  $k'$  which contains all the coefficients of $ F_i$  and is a finite Galois extension of $ k$.  Let  $G = G(L/k)$  be its Galois group and put $ m = \#G$.  Then  $G $ acts on  a polynomial ring  $L[X_1, \ldots ,X_n]$  such that $ X_i^g = X_i$  for all $ i$  and all $ g \in G$  that is,  $G$  acts on coefficients of an element in $ L[X_1, \ldots,  X_n]$.  Hence  
  $$  mX_i = \sum_{g \in G}X_i^g  = \sum_{g \in G}F_i^g(f_1^g, \ldots ,f_n^g) = \sum_{g \in G}F_i^g(f_1, \ldots ,f_n).$$   
Since $ \sum_{g \in G}F_i^g(Y_1,\ldots,Y_n) \in k[Y_1,\ldots,Y_n]$, it follows that $ \sum_{g \in G}F_i^g(f_1,\ldots,f_n)$ $\in k[f_1,\ldots,f_n]$.  Therefore $ X_i \in k[f_1,\ldots,f_n]$  because  $L$  has a characteristic zero.  So we may assume that $ k$  is algebraically closed.

(2)  Let  $k$  be a field, let  $k[X_1, \ldots, X_n]$ denote a polynomial ring  and let  $f_1, \ldots, f_n \in k[X_1, \ldots, X_n]$.  If  the Jacobian $\det\left(\dfrac{\partial f_i}{\partial X_i}\right) \in k^\times ( = k \setminus (0))$, then the $k[X_1, \ldots, X_n]$ is unramified over the subring $k[f_1, \ldots, f_n]$.
 Consequently $f_1, \ldots, f_n$  is algebraically independent over $k$.

 In fact, put  $T = k[X_1, \ldots, X_n]$  and $S = k[f_1, \ldots, f_n] ( \subseteq T)$.
We have an exact sequence by {\bf [9, (26.H)]} :
$$ \Omega_{S/k}\otimes_ST \mapright{v} \Omega_{T/k} \mapright{} \Omega_{T/S} \mapright{} 0,$$
where 
 $$v(df_i\otimes 1) = \sum_{j=1}^n\dfrac{\partial f_i}{\partial X_j}dX_j\ \ \ \ \ (1 \leq i \leq n).$$
  So $\displaystyle \det\left(\dfrac{\partial f_i}{\partial X_j}\right) \in k^\times$ implies that  $v$ is an isomorphism.  Thus  $\Omega_{T/S} = 0$  and hence  $T$  is unramified over  $S$ by {\bf [2, VI,(3.3)]} or {\bf [9]}.  Moreover  $K(T)$  is algebraic over  $K(S)$,  which means that  $f_1, \ldots, f_n$  are algebraically independent over  $k$.    
    \end{remark}

As a result of Theorem \ref{2.2. Theorem}, we have the following.

\begin{theorem}[The Jacobian Conjecture]\label{2.3. Theorem}  Let $ k$  be a field of characteristic zero, let  $k[X_1,\ldots,X_n]$  be a
 polynomial ring over  $k$, and let  $f_1, \ldots ,f_n$  be elements
 in  $k[X_1,\ldots,X_n]$.  Then  the Jacobian matrix $ (\partial f_i/\partial X_j)$  is invertible if and only if $ k[X_1,\ldots,X_n] = k[f_1,\ldots,f_n]$.
\end{theorem}

%%%%%%%%%%%%%%%
%%%%%%%%%%%%%%%%

\vspace{10mm}

\section{Generalization of The Jacobian Conjecture}

%%\markboth{\protect\footnotesize{}}{\protect\footnotesize{}} 

   The Jacobian Conjecture (Theorem \ref{2.3. Theorem}) can be generalized as follows.

\begin{theorem}\label{4.1. Theorem}   Let  $A$  be an integral
 domain whose quotient field  $K(A)$  is of characteristic
 zero.  Let $ f_1, \ldots ,f_n$  be elements of a
 polynomial ring $ A[X_1,\ldots,X_n]$  such that the
 Jacobian determinant $ det(\partial f_i/\partial X_j)$
  is a unit in $ A$.  Then
    $$A[X_1,\ldots,X_n] = A[f_1,\ldots,f_n].$$
\end{theorem}
 
\begin{proof}    It  suffices to prove $ X_1, \ldots ,X_n \in A[f_1, \ldots ,f_n]$.  We have  $ K(A)[X_1,\ldots,X_n] = K(A)[f_1,\ldots,f_n]$  by Theorem \ref{2.3. Theorem}.  Hence \
$$  X_1 = \sum c_{i_1 \cdots i_n}f_1^{i_1} \cdots f_n^{i_n}$$ 
with  $c_{i_1 \cdots i_n} \in K(A)$.  If we set
  $f_i = a_{i1}X_1 + \ldots + a_{in}X_n +$ $\mbox{(higher
 degree terms)}$, $ a_{ij} \in A$ ,  then the assumption
 implies that the determinant of a matrix  $(a_{ij})$  is a unit in $ A$.
 Let 
$$  Y_i = a_{i1}X_1 + \ldots + a_{in}X_n \ \ \ (1 \leq i \leq n).$$ \
Then  $A[X_1, \ldots ,X_n] = A[Y_1, \ldots ,Y_n]$
  and $ f_i = Y_i + \mbox{(higher degree terms)}$. 
 So to prove the assertion,  we can assume that without 
loss of generality the linear parts of $ f_1, \ldots ,f_n$  are
 $ X_1,\ldots,X_n$,  respectively.  Now we introduce
 a linear order in the set  $ \{(i_1,\ldots,i_n) \mid i_k \in {\bf Z}\}$  of lattice points in ${\bf R}^n$ (where ${\bf R}$ denotes the field of real numbers)  in the way :
 $(i_1, \ldots ,i_n) > (j_1, \ldots ,j_n)$  if (1) $ i_1 + \ldots + i_n > j_1 + \ldots + j_n $ or (2) $i_1 + \ldots + i_k > j_1 + \ldots + j_k$  and  $i_1 + \ldots + i_{k+1} = j_1 + \ldots + j_{k+1}, \ldots \ldots, i_1 + \ldots + i_n = j_1 + \ldots + j_n$.  We shall show that every $ c_{i_1 \ldots i_n}$  is in $ A$  by induction on the linear order
 just defined.  Assume that every  $ c_{j_1 \ldots j_n}$
  with  $(j_1, \ldots ,j_n) < (i_1, \ldots ,i_n)$  is
 in $ A$.  Then the coefficients of the polynomial \
$$  \sum c_{j_1 \cdots j_n}f_1^{j_1} \cdots f_n^{j_n}$$
are in $ A$, where the summation ranges over  $(j_1,\ldots,j_n) \geq  (i_1, \ldots ,i_n)$.  In this polynomial, the
 term  $X_1^{i_1} \cdots X_n^{i_n}$  appears once with
 the coefficient $ c_{i_1 \ldots i_n}$.  Hence  $ c_{i_1 \ldots i_n}$  must be an element of  $A$.  So  $X_1$  is in $ A[f_1, \ldots ,f_n]$.  Similarly $X_2, \ldots ,X_n$ 
 are in  $A[f_1,\ldots,f_n]$  and the assertion
 is proved completely. \end{proof}

\begin{corollary}\label{4.2. Corollary}{\rm (Keller's Problem)}    Let  $f_1, \ldots ,f_n$  be elements of a polynomial ring  ${\bf Z}[X_1,\ldots,X_n]$  over  ${\bf Z}$, the ring of integers.  If the Jacobian
 determinant $ det(\partial f_i/\partial X_j)$  is equal to either \ $\pm 1$, then  ${\bf Z}[X_1,\ldots,X_n] = {\bf Z}[f_1,\ldots,f_n]$.
\end{corollary}

   %%%%%
       
\bigskip
\bigskip
 \bigskip
\bigskip
\bigskip
\bigskip

\noindent
Department of Mathematics\\
 Faculty of Education\\
  Kochi University\\
   2-5-1 Akebono-cho, Kochi 780-8520\\
 JAPAN\\
 ssmoda@kochi-u.ac.jp

\end{document}